\documentclass[pdflatex,iicol]{sn-jnl}

\usepackage{graphicx}%
\usepackage{multirow}%
\usepackage{amsmath,amssymb,amsfonts}%
\usepackage{amsthm}%
\usepackage{mathrsfs}%
\usepackage[title]{appendix}%
\usepackage[dvipsnames]{xcolor}%
\usepackage{textcomp}%
\usepackage{manyfoot}%
\usepackage{booktabs}%
\usepackage{algorithm}%
\usepackage{algorithmicx}%
\usepackage{algpseudocode}%
\usepackage{listings}%

\usepackage{natbib}
\newcommand{\cmmnt}[1]{\ignorespaces}
\usepackage{upgreek} % For non-italic and "boldable" greek letters.
\usepackage{physics} % For derivatives.
\newcommand{\mvs}{max$(|v_\mathrm{solid}|)$ } % Shortcut.
\newcommand{\am}{$\alpha_\mathrm{max}$ } % Shortcut.
\usepackage{enumitem}
\usepackage{float}
\usepackage{relsize} % For \mathlarger command.

\usepackage{esint}
\usepackage{dirtytalk}
\usepackage[normalem]{ulem}
\usepackage{natbib}
\usepackage[font=scriptsize]{subfig}
\usepackage{epstopdf}
\usepackage{multirow}
\usepackage{mwe}
\usepackage{makecell}
\usepackage{hhline}
\usepackage{float}
\usepackage{bookmark}
\usepackage{tabu}
\usepackage{upgreek} % For non-italic and "boldable" greek letters.
\usepackage{relsize} % For \mathlarger command.
\usepackage{physics} % For derivatives.
\usepackage{amssymb} % For \square.
\usepackage{enumitem} % Change symbol for enumerate.
\makeatletter
\def\env@dmatrix{\hskip -\arraycolsep
	\let\@ifnextchar\new@ifnextchar
	\extrarowheight=0ex
	\array{*\c@MaxMatrixCols{>{\displaystyle}c}}}

\newenvironment{bdmatrix}
{\left[\env@dmatrix}
{\endmatrix\right]}
\newenvironment{Bdmatrix}
{\left\{\env@dmatrix}
{\endmatrix\right\}}
\makeatother

\theoremstyle{thmstyleone}%
\theoremstyle{thmstyletwo}%

\theoremstyle{thmstylethree}%

\begin{document}

\title[On the Calculation of the Brinkman Penalization Term in Density-Based Topology Optimization of Fluid-Dependent Problems]{On the Calculation of the Brinkman Penalization Term in Density-Based Topology Optimization of Fluid-Dependent Problems}

%%=============================================================%%
%% GivenName	-> \fnm{Joergen W.}
%% Particle	-> \spfx{van der} -> surname prefix
%% FamilyName	-> \sur{Ploeg}
%% Suffix	-> \sfx{IV}
%% \author*[1,2]{\fnm{Joergen W.} \spfx{van der} \sur{Ploeg} 
%%  \sfx{IV}}\email{iauthor@gmail.com}
%%=============================================================%%

\author[1]{\fnm{Mohamed} \sur{Abdelhamid}}%\email{maa.abdelhamid@gmail.com}

\author*[1]{\fnm{Aleksander} \sur{Czekanski}}\email{alex.czekanski@lassonde.yorku.ca}

\affil*[1]{\orgdiv{Department of Mechanical Engineering}, \orgname{York University}, \orgaddress{\street{4700 Keele St.}, \city{Toronto}, \postcode{M3J 1P3}, \state{Ontario}, \country{Canada}}}

%%==================================%%
%% Sample for unstructured abstract %%
%%==================================%%

\abstract{In topology optimization of fluid-dependent problems, there is a need to interpolate within the design domain between fluid and solid in a continuous fashion. In density-based methods, the concept of inverse permeability as of a volumetric force is utilized to enforce zero fluid velocity in non-fluid regions. This volumetric force consists of a scalar term multiplied by fluid velocity. This scalar term takes a value between two limits as determined by a convex interpolation function. The maximum inverse permeability limit is typically chosen through a trial and error analysis of the initial form of the optimization problem; such that the fields resolved resemble those obtained through an analysis of a pure fluid domain with a body-fitted mesh. In this work, we investigate the dependency of the maximum inverse permeability limit on mesh size and flow conditions through analyzing the Navier-Stokes equation in its strong and discretized finite element forms.}


%%\pacs[JEL Classification]{D8, H51}

%%\pacs[MSC Classification]{35A01, 65L10, 65L12, 65L20, 65L70}

\maketitle

\section{Introduction}
\label{sec:intro}

\subsection{Brinkman Penalization as a Design Parametrization Technique}
The first application of \textbf{t}opology \textbf{o}ptimization (TO) to fluid-dependent problems appeared in the seminal work of \cmmnt{Borrvall and Petersson} \cite{Borrvall2003}, where they addressed a pure fluid problem under Stokes flow conditions. Later, \cmmnt{Gersborg-Hansen et al.} \cite{Gersborg-Hansen2005} extended the work to Navier-Stokes equations. Although both works utilized the analogy of a 2D channel flow with varying thickness for design parametrization, the later recognized the similarity between this model and Brinkman equations of fluid flow in porous media \citep{Brinkman1947}. This similarity was also noted independently by \cmmnt{Evgrafov} \cite{Evgrafov2005} and \cmmnt{Guest and Pr\'{e}vost} \cite{Guest2006}, where the later directly used Darcy's law - a porous flow model - to introduce fluid flow in porous regions, hence freeing the topology optimization model from its two-dimensional channel assumption. In addition, the use of a porous flow model such as Darcy's law warranted a physical interpretation of porosity for intermediate densities. Consequently, this model could potentially be used to design porous media such as filters, and it's no longer \textit{a mere bridge} to interpolate between fluid and solid with the final goal of reaching only pure discrete designs \citep[p.~463]{Guest2006}. This now termed \textit{Brinkman penalization} is the de facto method for density-based topology optimization of fluid-dependent problems. In the remainder of this work, `Brinkman penalization' and `inverse permeability' are used interchangeably and our discussion is limited to finite element discretizations such that each finite element is parameterized using a single fluid design variable $\rho$.

Typically, Brinkman penalization is employed by appending a negative volumetric force to the body force and internal force terms in the Navier-Stokes momentum equation. This volumetric force is basically the Brinkman inverse permeability scalar function multiplied by the velocity vector, such that it has a scalar component in each spatial direction; i.e. $x$ and $y$ in 2D. This Brinkman penalization function is convex and ranges between a maximum and a minimum limit. It usually takes the following form which first appeared in \cite{Borrvall2003}:
\begin{equation}
	\alpha (\rho) = \alpha_\text{max} + \rho (\alpha_\text{min} - \alpha_\text{max}) \frac{1 + p_\alpha}{\rho + p_\alpha}
	\label{eq:alpha}
\end{equation} 

% \footnote{The authors recognize it's unorthodox to list an equation in the introduction section, however, we believe it's highly pertinent to the remaining literature survey which is focused on the details of this particular equation.}

\noindent where $\alpha_\text{max}$ and $\alpha_\text{min}$ are the maximum and minimum inverse permeability limits (also known as Brinkman penalization limits), $\rho$ is the fluid design variable ($\rho = 1$ for fluid and $\rho = 0$ for solid), and $p_\alpha$ is the Brinkman penalization interpolation parameter. This Brinkman penalization function is different from the somewhat \textit{analogous} \textbf{S}olid \textbf{I}sotropic \textbf{M}aterial with \textbf{P}enalization (SIMP) function used with topology optimization of solids in two aspects:
\begin{enumerate}[label=(\alph*)., labelindent=0pt, labelwidth=!]
\item Unlike SIMP, since the volumetric force term is appended to other non-zero terms, there is no fear of singularities - from a mathematical perspective - if the volumetric force term vanishes in pure fluid elements \citep[p.~469]{Guest2006}. Hence, $\alpha_\text{min}$ maybe taken as zero except when there is a physical need to be non-zero as in solving a two-dimensional problem with a finite out-of-plane thickness such as a microfluidic device, cf. \citep[p.~182]{Gersborg-Hansen2005}, \citep[p.~978]{Olesen2006}, and \citep[p.~5]{Alexandersen2023}.

\item A linear SIMP function in terms of the design variable means no penalization is imposed on intermediate elements, while in Brinkman penalization, a linear relation enforces severe penalization on intermediate elements.
\end{enumerate}

\subsection{Calculation of $\alpha_\text{max}$}

As for $\alpha_\text{max}$, it is typically selected just high enough to enforce near zero velocity in non-fluid elements, yet small enough not to introduce numerical instabilities. From early on, the significance of this maximum limit was recognized and its effect on the optimized design was discussed. In \citep[p.~102]{Borrvall2003}, the authors recognized the strong dependence of the objective function of power dissipation on this maximum limit, yet the optimized designs were found to be highly independent of that limit. In \citep[p.~184]{Gersborg-Hansen2005}, the authors chose the maximum limit so as to enforce a small flow rate in the solid regions in the range of two orders of magnitude lower than the maximum velocity in pure fluid regions.

In \cite{Guest2006}, the authors studied the effect of different magnitudes of the maximum limit by solving a sample problem through two models; a pure fluid model and their Darcy-Stokes model developed for TO. The resolved fluid velocity and pressure fields from the two models were then compared. They noted a linear relation between the permeability (i.e. $1/\alpha_\text{max}$ in this work) and the maximum velocity in the solid regions. A deterioration in the solution accuracy coincided with the loss of that linear relation, which occurred at low permeability values (i.e. equivalent to high $\alpha_\text{max}$ in our formulation). In our numerical experiments, we noticed a different behavior that is discussed in detail in Section \ref{sec:character}.

\cmmnt{Guest and Pr\'{e}vost} \cite{Guest2006} also hinted at the dependency of the Brinkman penalization limits on the mesh size utilized by calculating a certain permeability value as a function of mesh size. This value corresponded to equal diagonal terms between the Darcy and Stokes stiffness matrices, and was later used as an initial value for their implemented continuation technique. In contrast to \cite{Borrvall2003} and \cite{Gersborg-Hansen2005} which gradually raised $p_\alpha$ in Eq. \ref{eq:alpha} to introduce continuation, the authors in \cite{Guest2006} implemented what is analogous to $\alpha(\rho) = \rho \ \alpha_\text{max}$ directly and gradually raised $\alpha_\text{max}$ to introduce continuation.

In \cite{Olesen2006}, the authors calculated the proper maximum inverse permeability limit by looking at the streamlines in the resolved velocity field to estimate how much flow went through the solid structure, and also by looking at the relation of the objective function w.r.t. the maximum limit. They also mentioned that the maximum limit could be of use as a continuation tool in severe nonconvex problems similar to \cite{Guest2006}.

\cmmnt{Kreissl et al.} \cite{Kreissl2011} noted the independence of the maximum impermeability limit on the Reynolds number. They also noted the need for a relatively fine mesh for the pressure fields to match between the Brinkman-penalized Navier-Stokes and original Navier-Stokes with body-fitted mesh.

In recent literature on \textbf{t}opology \textbf{o}ptimization of \textbf{f}luid-\textbf{s}tructure \textbf{i}nteraction (TOFSI) problems, $\alpha_\text{max}$ is calculated by solving an analysis of an initial discrete design of the TO problem using a body-fitted mesh of segregated, non-overlapping fluid/solid domains. A parameter of interest obtained from this analysis is used as a benchmark against the same parameter calculated by analyzing the unified domain formulation with the Brinkman penalization term implemented. The maximum limit $\alpha_\text{max}$ is usually progressively increased by an order of magnitude until the two results match within a certain error margin, cf. \citep[p.~602]{Yoon2010} and \citep[p.~993]{Lundgaard2018}.

While the trial-and-error approach for selecting a proper $\alpha_\text{max}$ is acceptable for a single design problem, sometimes a need arises for calibrating $\alpha_\text{max}$ such that different mesh sizes and flow conditions produce the same behavior. In particular, there is usually a necessity for solving the TO problem using a relatively coarse mesh before committing to solving the final refined (hence costly) mesh, such as the need to calibrate some interpolation and projection parameters. In fact, the motivation for this study arose in the authors' work on density-based TOFSI problems; a multiphysics problem known for its tediously strong nonlinear and nonconvex behavior. We noticed that after calibrating some interpolation parameters on a relatively coarse mesh that is solvable within a reasonable time frame, the same parameters produced a different behavior when applied to the finer mesh needed for producing the final results.

\textbf{In this work}, we investigate the dependency of the Brinkman penalization term on the mesh size and the flow conditions. Through analyzing the Navier-Stokes equations in their PDE as well as discretized finite element forms, we propose proportionality relations to describe these dependencies. We solve a wide range of numerical experiments and use curve fitting to characterize these dependencies. The rest of this manuscript is organized as follows; in \textbf{Section \ref{sec:govern_eqns}}, we introduce the fluid flow governing equations and boundary conditions. In \textbf{Section \ref{sec:finite}}, we discuss the finite element discretization of the governing equations which provide valuable insights into the dependency of the Brinkman penalization maximum limit on the mesh size. In \textbf{Section \ref{sec:dependency}}, we analyze the fluid flow governing equations to deduce proportionality relations between the Brinkman penalization term and mesh size and flow conditions. In \textbf{Section \ref{sec:character}}, we use numerical experiments to verify and characterize the proportionality relations derived. Finally, in \textbf{Section \ref{sec:conc}}, we summarize our findings and present our concluding remarks.

\section{Governing Equations of Fluid Flow}
\label{sec:govern_eqns}

Before starting out investigation into the dependence of the Brinkman maximum limit on the mesh size and the flow conditions, we should first establish the governing equations of the problem at hand. Consider the Navier-Stokes equations in their incompressible, steady-state form \cite[p.~10]{Reddy2010}. The strong form of the PDEs modified for TO is as follows:
\begin{gather}
	\boldsymbol{\nabla} \cdot \mathbf{v} = 0, \label{eq:contin} \\[1em]
	\rho_f \, (\mathbf{v} \cdot \boldsymbol{\nabla}) \mathbf{v} = \boldsymbol{\nabla} \cdot \boldsymbol{\upsigma}^f + \mathbf{f}^f - \alpha(\rho) \mathbf{v} , \label{eq:ns} \\[1em] 
	\boldsymbol{\upsigma}^f = - p \mathbf{I} + \mu \left[ \boldsymbol{\nabla} \mathbf{v} + (\boldsymbol{\nabla} \mathbf{v})^T \right], \\[1em]
	\alpha(\rho) = \alpha_\text{max} + \rho \left( \alpha_\text{min} - \alpha_\text{max} \right) \frac{1 + p_\alpha}{ \rho + p_\alpha}. \label{eq:invrs_perm}
\end{gather}

\noindent where $\mathbf{v}$ is the fluid velocity, $\rho_f$ is the fluid density (a subscript $f$ is used to distinguish it from the design variables $\rho$), $\boldsymbol{\upsigma}^f$ is the Cauchy fluid stress tensor, $\mathbf{f}^f$ is the external fluid force (assumed zero and dropped in the remainder of this work), $p$ is the hydrostatic pressure, and $\mu$ is the fluid dynamic viscosity. The fluid momentum equation, Eq. \ref{eq:ns}, is appended with the Brinkman penalization term $- \alpha(\rho) \mathbf{v}$ as a volume force to enforce zero velocity in 0\% fluid elements while allowing for smooth interpolation between the artificial density limits; i.e. between solid and fluid. The Brinkman penalization interpolation parameter $p_\alpha$ is usually selected based on the physics of the problem at hand; e.g. Reynolds number in TOFSI problems \cite[p.~974]{Lundgaard2018}. A continuation scheme maybe used with $p_\alpha$ to avoid the optimizer getting stuck in local minima; cf. \citep[p.~96]{Borrvall2003} and \citep[p.~10]{Alexandersen2023}.

The essential boundary conditions are defined as follows:
\begin{alignat}{3}
	\text{Fluid No-slip:} & \qquad \mathbf{v} = 0 & \quad \text{on } & \mathit{\Gamma}_{\mathbf{v}_0}, \label{eq:bc_noslip} \\[1em]
	\text{Fluid Inlet:} & \qquad \mathbf{v} = \mathbf{v}_\text{in} & \quad \text{on } & \mathit{\Gamma}_{\mathbf{v}_\text{in}}, \label{eq:bc_v_in} \\[1em]
	\text{Fluid Outlet:} & \qquad p = 0 & \quad \text{on } & \mathit{\Gamma}_{\mathbf{v}_\text{out}}. \label{eq:bc_p_out}
\end{alignat}

\noindent where $\mathbf{v}_\text{in}$ is the prescribed inlet velocity at the inlet boundary $\mathit{\Gamma}_{\mathbf{v}_\text{in}}$, and $\mathit{\Gamma}_{\mathbf{v}_\text{out}}$ is the outlet boundary with a prescribed zero pressure applied. Note that the fluid no-slip boundary condition in Eq. \ref{eq:bc_noslip} is only defined on the remaining external domain boundaries $\mathit{\Gamma}_{\mathbf{v}_0}$ aside from the inlet and outlet boundaries such that  $ \mathit{\partial \Omega}_f = \mathit{\Gamma}_{\mathbf{v}_0} \cup \mathit{\Gamma}_{\mathbf{v}_\text{in}} \cup \mathit{\Gamma}_{\mathbf{v}_\text{out}} $. The volume force term appended to the fluid momentum in Eq. \ref{eq:ns} automatically enforces a no-slip condition wherever needed within the solid domain and its boundaries.

\section{Finite Element Formulations}
\label{sec:finite}

To study the dependence of the Brinkman penalization upper limit $\alpha_\text{max}$ on the mesh size, we must take a closer look at the discretized weak form of the Navier-Stokes and continuity equations. We implement the \textit{standard Galerkin method of weighted residuals} where the test/weight functions are the same as the interpolation/shape functions. The resulting model is of the \textit{velocity-pressure} (or \textit{mixed}) type where both velocity and pressure are solved for simultaneously. To satisfy the \textit{Ladyzhenskaya-Babuska-Brezzi} condition, cf. \cite[p.~176]{Reddy2010}, \textit{P2P1 Lagrangian} finite elements (i.e. 9 velocity nodes and 4 pressure nodes) are used with a low to moderate Reynolds number to avoid using stabilization techniques that artificially - but not necessarily accurately - dampen the discontinuities. We employ regular rectangular meshing of equal size, the mesh size $h$ denotes the length of the finite element side. The continuity equation (Eq. \ref{eq:contin}) is typically weighted by the pressure shape function $\mathbf{\Phi}$ while the momentum equation (Eq. \ref{eq:ns}) is typically weighted by the velocity shape function $\mathbf{\Psi}$. The resulting finite element system in 2D \textit{on the elemental level} is as follows:

\begin{gather}
	\underbrace{\begin{bmatrix}
			2 \mathbf{K}_{11} + \mathbf{K}_{22} + \mathbf{C(v)} & \mathbf{K}_{12} & - \mathbf{Q}_1 \\[4pt]
			\mathbf{K}_{21} & \mathbf{K}_{11} + 2 \mathbf{K}_{22} + \mathbf{C(v)} & - \mathbf{Q}_2 \\[4pt]
			- \mathbf{Q}_1^T & - \mathbf{Q}_2^T & \mathbf{0}
	\end{bmatrix}}_{\text{Conservation of Momentum and Mass}} \nonumber \\[1em]
	\begin{Bmatrix}
		\mathbf{\hat{v}}_1 \\[4pt]
		\mathbf{\hat{v}}_2 \\[4pt]
		\mathbf{\hat{p}}
	\end{Bmatrix} +
	%%%
	\underbrace{\begin{bmatrix}
			\mathbf{A} & \mathbf{0} & \mathbf{0} \\[4pt]
			\mathbf{0} & \mathbf{A} & \mathbf{0} \\[4pt]
			\mathbf{0} & \mathbf{0} & \mathbf{0}
	\end{bmatrix}}_{\substack{\text{Brinkman} \\ \text{Penalization}}}
	\begin{Bmatrix}
		\mathbf{\hat{v}}_1 \\[4pt]
		\mathbf{\hat{v}}_2 \\[4pt]
		\mathbf{\hat{p}}
	\end{Bmatrix} = \mathbf{0}.
	\label{eq:main_fea}
\end{gather}

The coefficient matrices in the finite element form are defined as ($\int_{-1}^{+1} \int_{-1}^{+1}$ and $\dd{\xi} \dd{\eta}$ are implied)\footnote{Summation is implied on repeated indices in $\mathbf{C(v)}$ but not in $\mathbf{K}_{ij}$.}:
\begingroup
\allowdisplaybreaks
\begin{gather}
	\mathbf{K}_{ij}  = \mu \, \pdv{\mathbf{\Psi}}{x_i} \pdv{\mathbf{\Psi}}{x_j} ^T |\mathbf{J}|, \label{eq:2} \\[1em]
	%%%
	\mathbf{C(v)}  = \rho \, \mathbf{\Psi} \left[ \left( \mathbf{\Psi}^T \mathbf{\hat{v}}_i \right) \pdv{\mathbf{\Psi}}{x_i} ^T \right] |\mathbf{J}|, \\[1em]
	%%%
	\mathbf{Q}_i  = \pdv{\mathbf{\Psi}}{x_i} \, \mathbf{\Phi}^T |\mathbf{J}|, \\[1em]
	%%%
	\mathbf{A}  = \alpha(\rho) \mathbf{\Psi} \, \mathbf{\Psi}^T |\mathbf{J}|.
	\label{eq:brnk_fea}
\end{gather} 
\endgroup

\noindent where $\mathbf{\hat{v}}_1$ and $\mathbf{\hat{v}}_2$ are the nodal velocities in $x$ and $y$, respectively, and $\mathbf{\hat{p}}$ is the nodal pressures. $|\mathbf{J}|$ is the Jacobian determinant and $\xi$ and $\eta$ are the natural coordinates. No externally applied nodal fluid forces are used in this work as fluidic boundary conditions (Eqs. \ref{eq:bc_noslip} to \ref{eq:bc_p_out}) are implemented directly by setting nodal velocities/pressures to their appropriate values (i.e. strong, point-wise enforcement). Appropriate global assembly of Eq. \ref{eq:main_fea} is implemented and the resulting nonlinear system is solved using the undamped Newton-Raphson method \cite[p.~190]{Reddy2010}.

In the next section, we analytically investigate the dependence (or independence) of \am on the mesh size and the flow conditions.

\section{Analytical Derivation of the Dependence of $\alpha_\text{max}$ on Mesh Size and Flow Conditions}
\label{sec:dependency}

In order to establish the basis of our investigation, we take a look at two sets of parameters; namely the maximum state variable errors in 100\% fluid regions, and the maximum velocity in pure solid regions. We consider two perspectives; \textbf{(i)} the suitability of these parameters in measuring how the Brinkman penalized Navier-Stokes approximates the pure Navier-Stokes, and \textbf{(ii)} the easiness of investigating either set of parameters. In the following, we state our argument for the validity of using either set of parameters:
\begin{enumerate}[labelindent=0pt, labelwidth=!]
	\item The errors in the velocity and pressure fields resolved in 100\% fluid regions in comparison to those resolved using a pure fluid model. One of the main indicators of the validity of the Brinkman penalization model is that it should produce similar state fields to what is produced from a pure fluid model. This is even more critical in multiphysics problems whose behavior depends on the state variables such as structural compliance in TOFSI.
	
	\item The velocity in the solid regions is a good and direct indication of the validity of the Brinkman penalization model in simulating porous flow in solid media. In fact, it was one of the early parameters used in calibrating $\alpha_\text{max}$ as in \cite{Guest2006}.	
\end{enumerate}

Now that we established the validity of choosing either option from a representation point of view, next we look at the complexity of using either option from a mathematical equation-based perspective. The \textbf{first option} is a bit tricky to utilize as $\alpha_\text{max}$ does not have a direct influence on the 100\% fluid regions, instead the errors in the fluid state variables are reduced by minimizing the flow velocity in the pure solid regions hence directing the entire flow to the pure fluid regions and increasing the similarity to the results of a pure fluid model. Notice that in these discussions, we are looking at a special case, that is the existence of only discrete densities; either $\rho = 1$ in 100\% fluid regions or $\rho = 0$ in 100\% solid regions. The \textbf{second option}, i.e. velocity in solid regions, can be easily deduced through looking at the diagonal terms in Eq. \ref{eq:main_fea}.

In the following subsections, we discuss the dependency of the maximum inverse permeability limit on mesh size and flow conditions using the maximum velocity in solid regions\footnote{Typically, max$(|v_\mathrm{solid}|)$ should be scaled w.r.t. a nominal velocity characterizing the flow such as the characteristic $v_c$. In this work, $v_c$ is only changed by less than an order of magnitude, hence this scaling is not discussed further.} - designated as max$(|v_\mathrm{solid}|)$ - as a criterion through looking at the diagonal terms in the discretized finite element equations.

\subsection{Dependence of $\alpha_\text{max}$ on Mesh Size}
\label{ssec:alpha_h_anltc}

A closer look at the diagonal matrices in Eq. \ref{eq:main_fea} reveals the dependency of the original Navier-Stokes terms (i.e. $\mathbf{K}_{ij}$ and $\mathbf{C(v)}$) on the mesh size $h$ through the derivatives of the shape functions w.r.t. the global coordinates $\pdv*{\mathbf{\Psi}}{x_i}$. Recall that these derivatives are obtained as follows:
\begin{gather}
	\begin{Bdmatrix}
		\pdv{\mathbf{\Psi}}{x} \\[10pt] \pdv{\mathbf{\Psi}}{y}
	\end{Bdmatrix} = \left[ \mathbf{J} \right]^{-1} \begin{Bdmatrix}
		\pdv{\mathbf{\Psi}}{r} \\[10pt] \pdv{\mathbf{\Psi}}{s}
	\end{Bdmatrix}, \\[1em]
	\left[ \mathbf{J} \right]^{-1} = \frac{1}{|\mathbf{J}|} \begin{bdmatrix}
		+J_{2,2} & -J_{1,2} \\[8pt] -J_{2,1} & +J_{1,1}
	\end{bdmatrix}.
\end{gather}

\noindent where $J_{i,j}$ are the original Jacobian matrix components. Notice that, unlike the derivatives of $\mathbf{\Psi}$ w.r.t. the natural coordinates $\xi$ and $\eta$, $\pdv*{\mathbf{\Psi}}{x_i}$ is dependent on the mesh size $h$ through the components of the Jacobian matrix (in the numerator) and through the Jacobian determinant (in the denominator). This dependency can be characterized in closed form for the special case of regular, square meshing. The Jacobian matrix is calculated as follows:
\begin{equation}
	\left[ \mathbf{J} \right] = \begin{bdmatrix} 
		\pdv{\mathbf{\Psi}}{\xi} ^T \mathbf{\hat{x}} & \quad 		\pdv{\mathbf{\Psi}}{\xi} ^T \mathbf{\hat{y}} \\
		\pdv{\mathbf{\Psi}}{\eta} ^T \mathbf{\hat{x}} & \quad 		\pdv{\mathbf{\Psi}}{\eta} ^T \mathbf{\hat{y}}
	\end{bdmatrix}
\end{equation}

\noindent where the elemental nodal coordinates $\mathbf{\hat{x}}$ and $\mathbf{\hat{y}}$ are linearly proportional to the mesh size $h$ for the special case of regular, square finite elements.
In addition, the Jacobian determinant is typically related to the finite element's area (i.e. related to $h^2$). For a square element, the Jacobian determinant is known to be one fourth the element's area when evaluated anywhere within the element. This means that every $\pdv*{\mathbf{\Psi}}{x_i}$ is linearly proportional to the reciprocal of the mesh size, i.e. $1/h$. Again, the strength and regularity of this dependency depends on how distorted the element is from the ideal square shape.

On the other hand, the Brinkman penalization contribution to Eq. \ref{eq:main_fea} - namely matrix $\mathbf{A}$ - is independent of this parameter as it lacks any terms containing $\pdv*{\mathbf{\Psi}}{x_i}$. Hence, while the original Navier-Stokes terms change with different mesh sizes, the Brinkman penalization term does not.

From Eqs. \ref{eq:main_fea}-\ref{eq:brnk_fea}, it can be noted that the inverse permeability $\alpha(\rho)$ should be inversely proportional to $h^2$ (through $\mathbf{K}_{i,j}$ which contains two $\pdv*{\mathbf{\Psi}}{x_i}$ derivatives) and inversely proportional to $h$ (through $\mathbf{C(v)}$ which contains one $\pdv*{\mathbf{\Psi}}{x_i}$ derivative). Hence, the following relation between $\alpha_\text{max}$ and $h$:
\begin{equation}
	\alpha_\text{max} \propto \frac{1}{h^2} \qquad \& \qquad   \alpha_\text{max} \propto \frac{1}{h}
\end{equation}

\subsection{Dependence of $\alpha_\text{max}$ on Flow Conditions}
\label{ssec:kmax_re_anltc}

The dependence of the Brinkman penalization maximum limit $\alpha_\text{max}$ on the Reynolds number $Re$ can be investigated through looking at the non-dimensionalized form of Navier-Stokes equations. Following the treatment by \cite[p.~430]{leal2007advanced}, it is possible to non-dimensionalize Navier-Stokes equations w.r.t. the Reynolds number when under the assumptions of incompressible fluid, steady-state, and negligible body forces. Consider the following relations:
\begingroup
\allowdisplaybreaks
\begin{gather}
	\mathbf{v^*} = \frac{\mathbf{v}}{v_c}, \label{eq:nondmnsl_v} \\[1em]
	\boldsymbol{\nabla^*} = L_c \boldsymbol{\nabla}, \\[1em]
	p^* = \frac{p}{\rho_f \ {v_c}^2}, \\[1em]
	Re = \frac{v_c \ L_c \ \rho_f}{\mu}, \\[1em]
	\alpha^*(\rho) = \frac{\alpha(\rho) \ {L_c}^2}{\mu},
	\label{eq:nondmnsl_alpha} \\[1em]
	\alpha^* (\rho) = \alpha^*_\text{max} \left( 1 - \rho \frac{1 + p_\alpha}{\rho + p_\alpha} \right).
\end{gather} 
\endgroup

\noindent where the dimensionless form of each variable is designated with an asterisk superscript as in $\square^*$. $v_c$ is a characteristic velocity, taken in this work as the maximum inlet velocity in a parabolic laminar profile. $L_c$ is the characteristic length, taken as the width of the inlet boundary $\mathit{\Gamma}_{\mathbf{v}_\mathrm{in}}$. The relation in Eq. \ref{eq:nondmnsl_alpha} has been mentioned in relevant literature in some form, cf. \citep[p.~183]{Gersborg-Hansen2005}, \citep[p.~978]{Olesen2006}, and \citep[p.~598]{Yoon2010}. In that sense, Darcy's number $Da$ is equivalent to $1/\alpha^*_\text{max}$, both dimensionless. Generally, Darcy's number is related to a characteristic length that is relevant to the porous medium microstructure. It will be shown later in Section \ref{sec:character} that the characteristic length in Eq. \ref{eq:nondmnsl_alpha} should be different from, yet somehow related to, $L_c$.

By implementing Eqs. \ref{eq:nondmnsl_v}-\ref{eq:nondmnsl_alpha}, the Brinkman-penalized Navier-Stokes equations are non-dimensionalized in the following form:
\begin{equation}
	\begin{split}
		\frac{\rho_f \ {v_c}^2}{L_c} \, (\mathbf{v^*} \cdot \boldsymbol{\nabla^*}) \mathbf{v^*} =  - \frac{\rho_f \ {v_c}^2}{L_c} \, \boldsymbol{\nabla^*} \ p^* 	+ \\[1em]
		\frac{\mu \ v_c}{{L_c}^2} {\nabla^*}^2 \ \mathbf{v^*} - \frac{\mu \ v_c}{{L_c}^2} \alpha^*(\rho) \ \mathbf{v^*}
	\end{split}
\end{equation}

\noindent which can be rearranged through a multiplication by $L_c / \rho_f \ {v_c}^2$ as follows:
\begin{equation}
	\begin{split}
		(\mathbf{v^*} \cdot \boldsymbol{\nabla^*}) \mathbf{v^*} = - \boldsymbol{\nabla^*} \ p^* + \frac{1}{Re} {\nabla^*}^2 \ \mathbf{v^*} \\[1em]
		-   \frac{1}{Re} \alpha^*(\rho) \ \mathbf{v^*}
		\label{eq:nondmnsl_ns}
	\end{split}
\end{equation}

Similarly to the discussion in Section \ref{ssec:alpha_h_anltc}, we look at the diagonal terms in the finite element form of Eq. \ref{eq:nondmnsl_ns}. It appears that it's difficult to completely isolate $Re$ and its components, i.e. $v_c, \ \mu, \ \mathrm{and} \ L_c$, in a single term. Hence, although it might appear that $\alpha^*(\rho)$, hence $\alpha^*_\text{max}$, is independent of $Re$, $\alpha_\text{max}$ has the following dependencies:
\begin{equation}
	\alpha_\text{max} \propto \mu \quad \& \quad \alpha_\text{max} \propto \frac{1}{{L_c}^2} \quad \& \quad \alpha_\text{max} \propto v_c
	\label{eq:kmax_mu_lc_vc}
\end{equation}

In Eq. \ref{eq:kmax_mu_lc_vc}, the first two relations come from Eq. \ref{eq:nondmnsl_alpha} while the third one comes from the existence of velocity components in the convective term on the left hand side of Eq. \ref{eq:nondmnsl_ns} ($\mathbf{C(v)}$ in the finite element form). In other words, $\alpha_\text{max}$ is independent of $\rho_f$ but dependent on $v_c, \ \mu, \ \mathrm{and} \ L_c$.

In the next section, with the aid of numerical experiments, we focus on verifying the validity of the derived dependencies and on calculating the numerical values of the coefficients of proportionality derived earlier. Note that these coefficients of proportionality are only valid for the design problem discussed in this work. Nonetheless, we show that only a small number of data points is needed to calculate these coefficients for other problems.

\section{Characterizing the Dependency of $\alpha_\text{max}$ on Mesh Size and Flow Conditions}
\label{sec:character}

In this section, through numerical experiments, we aim to prove the validity of the derived dependencies in Section \ref{sec:dependency}; namely the dependency of the Brinkman penalization maximum limit $\alpha_\text{max}$ on the mesh size $h$, the fluid dynamic viscosity $\mu$, the characteristic length $L_c$, and the characteristic velocity $v_c$ as well as its independency of the fluid density $\rho_f$. In addition, we calculate exact numerical relations that describe these dependencies through curve fitting.

\begin{figure}[b]
	\centering
	\includegraphics[width=0.5\textwidth]{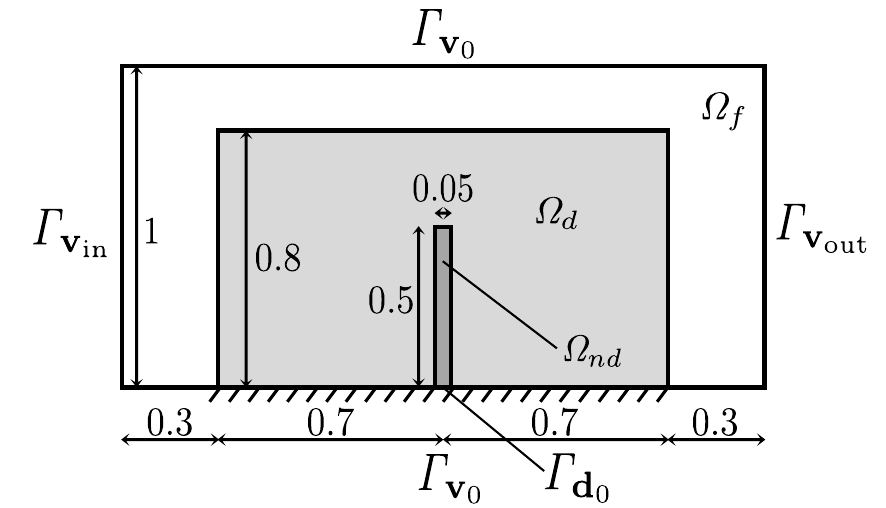}
	\caption{The \textit{modified} beam in a channel design problem as described in \cite{Lundgaard2018}.}
	\label{fig:mdfd_bm_in_a_chnl}
\end{figure}

To generate the data we use for curve fitting, we solve the governing equations of the Navier-Stokes equations equipped with the Brinkman penalization term. The problem to be solved is an initial design of the \textit{modified} beam in a channel problem described in Fig. \ref{fig:mdfd_bm_in_a_chnl}. The \textit{original} version of this problem was first discussed in a TOFSI context in \ \cite[p.~610]{Yoon2010} and has been used later as a benchmark problem in a number of works on TOFSI. It was later modified in \ \cite{Lundgaard2018}, hence the \textit{modified} designation, to increase the relative size of the design domain to the whole computational domain, rescale it from the micro to the macro scale, and generally strengthen the fluid-structure dependency.

As detailed in Fig. \ref{fig:mdfd_bm_in_a_chnl}, the problem features a 0.8 x 1.4 m rectangular design space (light gray) placed inside a 1 x 2 m rectangular channel. To avoid trivial solutions, a 0.05 x 0.5 m non-design solid beam (dark gray) is placed within the design space to force the optimizer to reach a more sophisticated solution than a simple bump at the bottom of the channel. The problem is solved for an initial discrete design such that $\rho = 0$ in $\Omega_d$ and $\Omega_{nd}$ and $\rho = 1$ in $\Omega_f \backslash \{\Omega_d \cup \Omega_{nd} \}$. Recall that, in this work, $\rho$ is defined as a fluid, not a solid, design variable.

The top and bottom surfaces of the channel $\mathit{\Gamma}_{\mathbf{v}_0}$ have a no-slip condition applied. A fully-developed, parabolic laminar flow profile of a maximum velocity $v_c$ is applied at the inlet $\mathit{\Gamma}_{\mathbf{v}_\text{in}}$ on the left, and a zero pressure condition is applied at the outlet $\mathit{\Gamma}_{\mathbf{v}_\text{out}}$ on the right. The bottom surface of the design and non-design spaces $\mathit{\Gamma}_{\mathbf{d}_0}$ is fixed to a ground structure. Note that even though this is a TOFSI problem, we are only concerned with the fluid flow analysis in this discussion.

As discussed earlier, the characteristic length $L_c$ is taken as the width of the entry boundary $\mathit{\Gamma}_{\mathbf{v}_\mathrm{in}}$ on the left. All the dimensions shown in Fig. \ref{fig:mdfd_bm_in_a_chnl} are scaled linearly with $L_c$. Unless otherwise noted, these default values are used:
\begingroup
\allowdisplaybreaks
\begin{gather}
	v_c = 1 \ \mathrm{m / s}, \\[1em]
	\rho_f = 1 \ \mathrm{kg / m^3}, \\[1em]
	\mu = 1 \ \mathrm{Pa \cdot s},	\\[1em]
	L_c = 1 \ \mathrm{m}, \\[1em]
	h = 0.01 \ \mathrm{m}, \\[1em]
	\alpha_\text{min} = 0 \ \mathrm{kg/m \cdot s}.
\end{gather}
\endgroup

First, we need to look at the maximum velocity in solid regions and the maximum state variable errors in fluid regions for a range of $\alpha_\text{max}$. Since such a study is a fluid flow analysis, it could be performed in a commercial software such as COMSOL Multiphysics \citep{comsol6} by employing the ``parametric sweep" feature. The Brinkman penalization term is easily implemented within the laminar flow model by using a volumetric force node and adding the scalar terms $- \alpha_\text{max} \ u$ and $- \alpha_\text{max} \ v$ in the $x$ and $y$ directions, respectively, where $u$ and $v$ are the $x$ and $y$ velocities as defined in the software.

\begin{figure*}[th!]
	\subfloat{\includegraphics[width=0.5\textwidth]{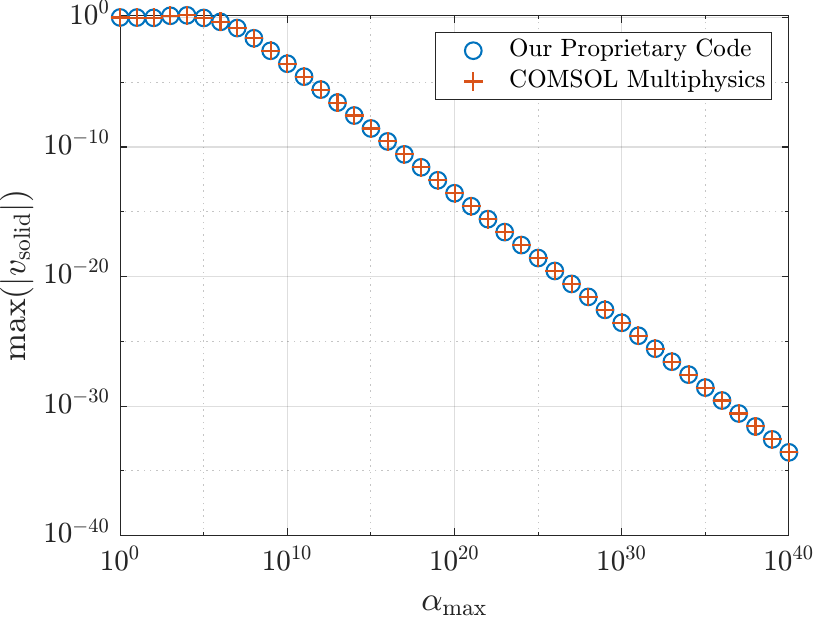}}
	\subfloat{\includegraphics[width=0.5\textwidth]{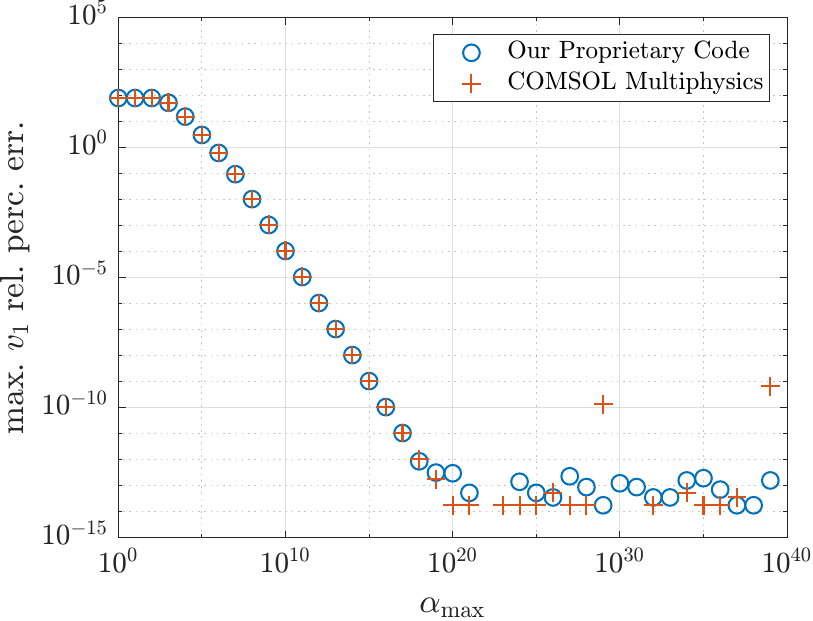}} \\
	\subfloat{\includegraphics[width=0.5\textwidth]{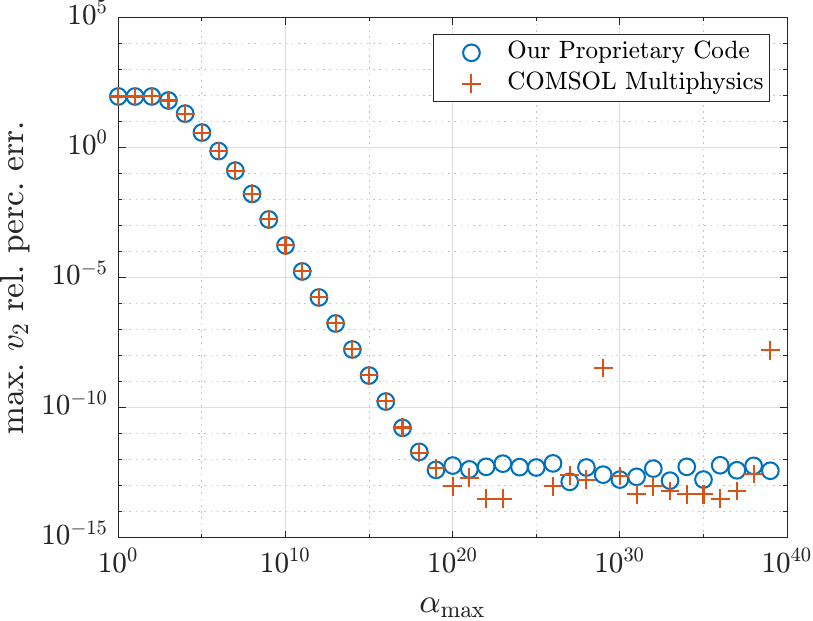}}
	\subfloat{\includegraphics[width=0.5\textwidth]{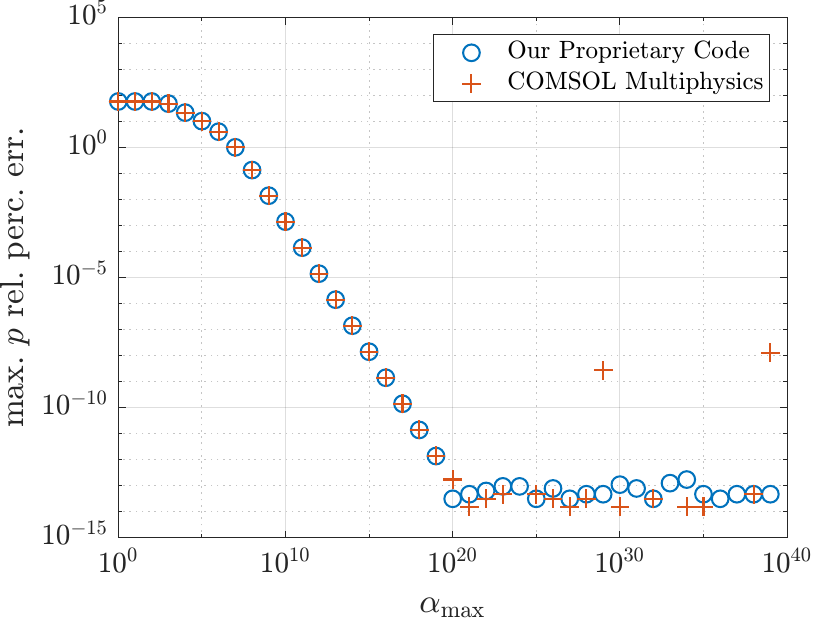}}
	\caption{Effect of $\alpha_\text{max}$ on the state variables in the solid and fluid domains. Results are obtained from our proprietary code in MATLAB as well as the commercial software COMSOL Multiphysics.}
	\label{fig:mtlb_vs_cmsl}
\end{figure*}

We solved the same problem using COMSOL as well as a proprietary code we wrote in MATLAB and the results are shown in Fig. \ref{fig:mtlb_vs_cmsl}. A mesh size of $h = 0.005$ m is utilized, resulting in a total of 80,000 finite elements. The problem is solved for the following range of $\alpha_\text{max}$ values; 0, $1e1$, $1e2$, ..., $1e40$.

The results obtained display a linear log-log relation between the Brinkman penalization maximum limit $\alpha_\text{max}$ and the maximum velocity in the 100\% solid regions. In contrast to the work of \ \citep[p.~471]{Guest2006}, we noted two differences; \textbf{(i)} we experience a linear relation between the log of the values, not the values themselves, and \textbf{(ii)} we don't notice any loss of accuracy at the high end of $\alpha_\text{max}$ even at considerably high values. We conjecture the reason for the first discrepancy is that \cite{Guest2006} addressed Stokes flow which neglected the nonlinear convection term while we are solving full Navier-Stokes equations. As for the second discrepancy, we conjecture this deterioration in accuracy is related to one or a combination of the following reasons; \textbf{(a)} the use of iterative solvers for the governing equations without proper preconditioners and tight convergence criteria, and \textbf{(b)} the use of stabilization techniques not calibrated for their newly-developed Darcy-Stokes model (cf. \ \citep[p.~1238]{Kreissl2011}). On the other hand, loss of linearity occurs understandably at low $\alpha_\text{max}$ values as the Brinkman-penalized Navier-Stokes model loses its accurate representation of the impermeable solid domain.

We also notice that the maximum absolute relative percentage errors in the pure fluid state variables maintain a linear log-log relation with $\alpha_\text{max}$ up till a certain limit ($\alpha_\text{max} \approx 1e18$ that is equivalent to max$(|v_\mathrm{solid}|) \approx 1e-12$), after which the values plateau at an almost constant level. We could argue that beyond this limit, no benefit is gained from using a higher $\alpha_\text{max}$ value. Hence, for the following results, we only run each study up till a maximum value of $\alpha_\text{max} = 1e20$.

In the following subsections, we discuss the dependency of $\alpha_\text{max}$ on each parameter individually.

\subsection{Relation between $\alpha_\text{max}$ and $h$}
\label{ssec:kmax_h_chrctr}

\begin{figure*}[t!]
	\centering
	\subfloat[Original.]{\includegraphics[width=0.5\textwidth]{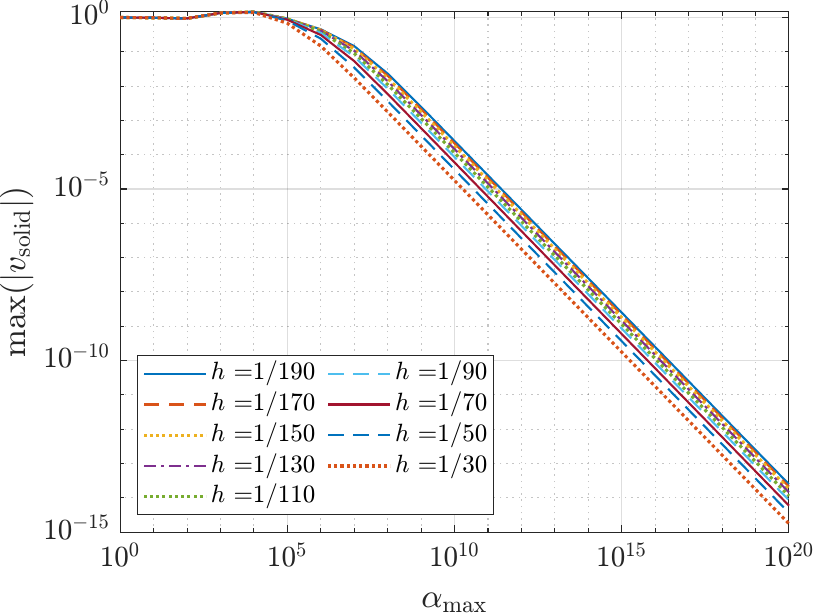}}
	\subfloat[Zoomed in.]{\includegraphics[width=0.5\textwidth]{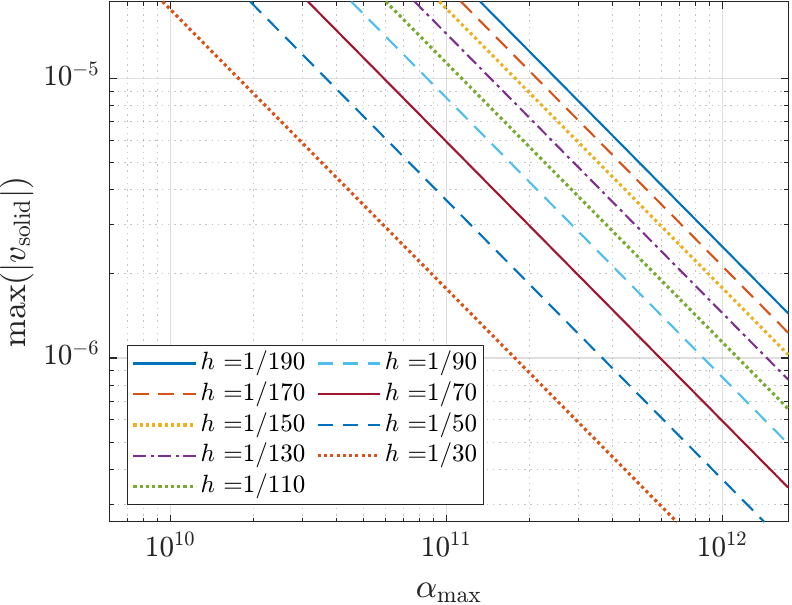}}
	\caption{Maximum velocity in the solid regions vs $\alpha_\text{max}$ for different $h$ values on a log-log scale.}
	\label{fig:kmax_vs_h} 
\end{figure*}

The first set of results concerns the dependency of $\alpha_\text{max}$ on the mesh size $h$. A fluid flow analysis is run for all combinations of the following values; $h = 1/30, 1/50, 1/70, ..., 1/190$ and $\alpha_\text{max} = 0, 1e0, 1e1,$ $..., 1e20$. The extracted results are presented in Fig. \ref{fig:kmax_vs_h}. It can be noted that for max$(|v_\mathrm{solid}|) \approx 1e-2$ and smaller, the log-log relation is linear.

\begin{figure*}[t!]
	\captionsetup{margin={10pt}}
	\centering
	\subfloat[A log-log scale is used. From bottom up, max$(|v_\mathrm{solid}|) = 1e-2, 1e-3, ..., 1e-12$.]{\includegraphics[width=0.5\textwidth]{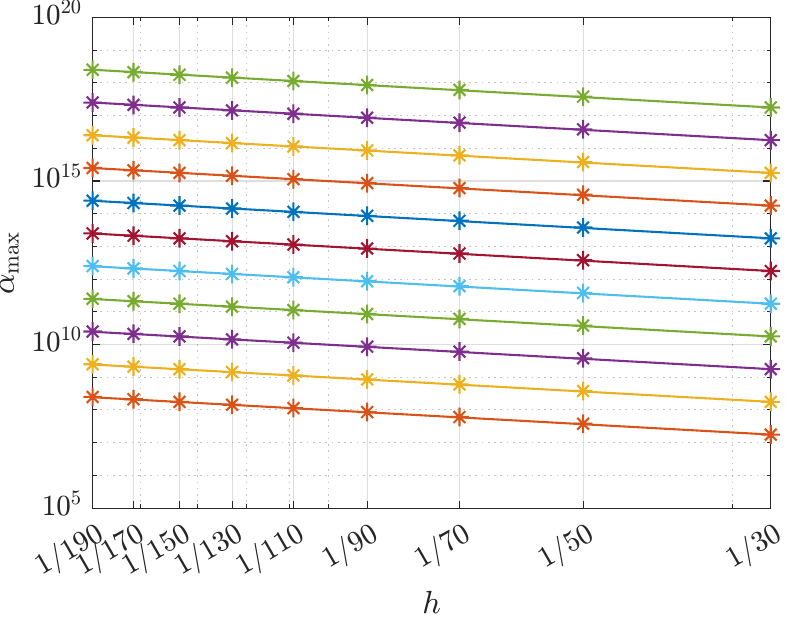}}
	\subfloat[A single case of max$(|v_\mathrm{solid}|) = 1e-12$ with error bars calculated as abs($\alpha_\text{max}|$\textsubscript{Data} - $\alpha_\text{max}|$\textsubscript{Eq.}).] {\includegraphics[width=0.5\textwidth]{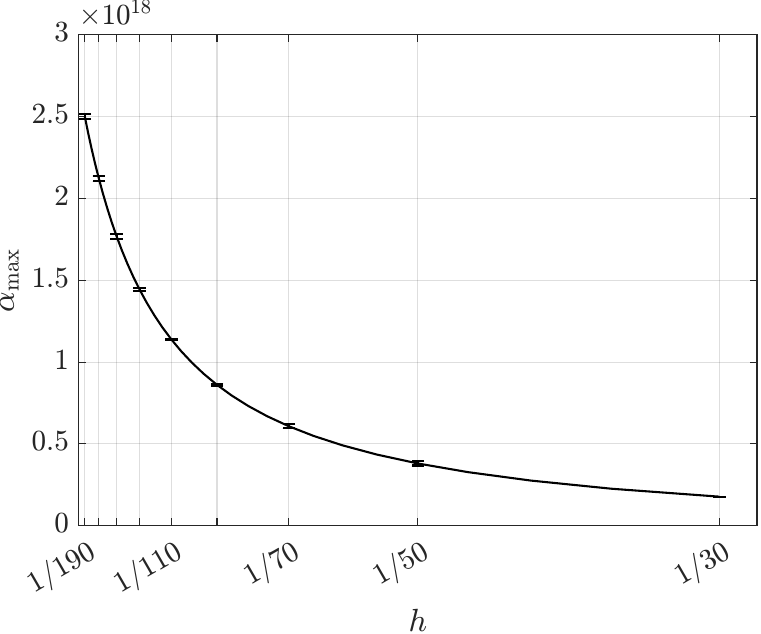}}
	\caption{Comparison of Eq. \ref{eq:kmax_vs_h_crv_ftng} (asterisks) to data points from the numerical experiments (solid lines).}
	\label{fig:kmax_vs_h_crv_ftng_chck}
\end{figure*}

In order to characterize the relation between $\alpha_\text{max}$ and $h$, we use curve fitting in order to calculate an expression for $\alpha_\text{max}$ as a function of $h$ and  max$(|v_\mathrm{solid}|)$. For curve fitting, we limit our data to only 6 points (3 points along h and 2 points along max$(|v_\mathrm{solid}|)$); which are all combinations of $h = 1/30, 1/110, 1/190$ and $\alpha_\text{max} = 1e8, 1e20$. We emphasize that for the curve fitting to be accurate, it is better for the data points used for fitting to be spanning the range of interest for each parameter. We note also that this choice of $\alpha_\text{max}$ ensures that \mvs $\le 1e-2$, hence within the linear portion of the log-log relation as presented in Fig. \ref{fig:kmax_vs_h}. Using curve fitting, the following relation is obtained:
\begin{equation}
	\alpha_\text{max} = 10^{-q} \left( \frac{31.32}{h^2} + \frac{7635}{h} - 8.039e04 \right)
	\label{eq:kmax_vs_h_crv_ftng}
\end{equation}

\noindent where $q$ is the exponent of the intended maximum velocity in the solid regions; i.e. max$(|v_\mathrm{solid}|) = 10^q$. To check the soundness of this relation, we compare Eq. \ref{eq:kmax_vs_h_crv_ftng} to the original set of data points for $h = 1/30, 1/50, 1/70, ..., 1/190$ and $\alpha_\text{max} = 1e8, 1e9, 1e10, $ $..., 1e20$. The comparison is presented in Fig. \ref{fig:kmax_vs_h_crv_ftng_chck} showing good agreement with a maximum error of 3.4\% for the case of max$(|v_\mathrm{solid}|) = 10^{-12}$.

\begin{figure*}[th!]
	\centering
	\subfloat[Original.]{\includegraphics[width=0.5\textwidth]{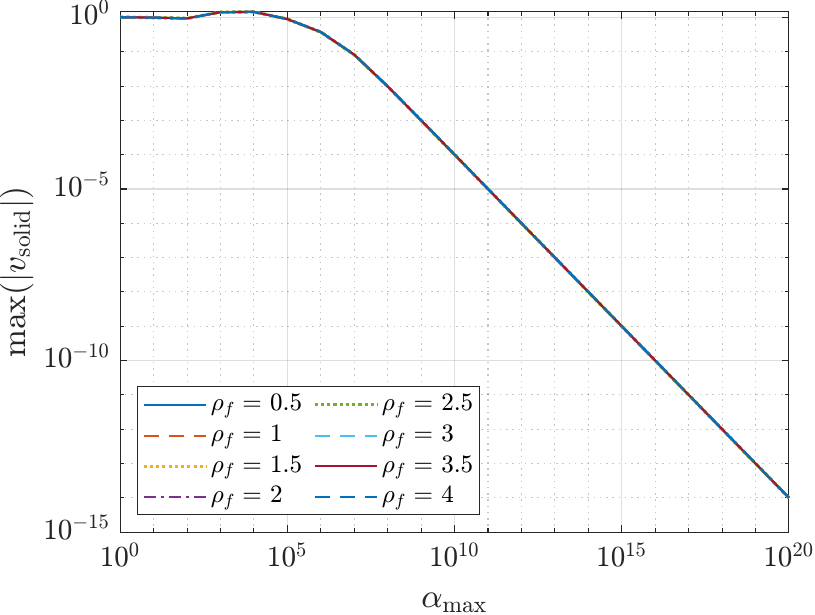}}
	\subfloat[Zoomed in.]{\includegraphics[width=0.5\textwidth]{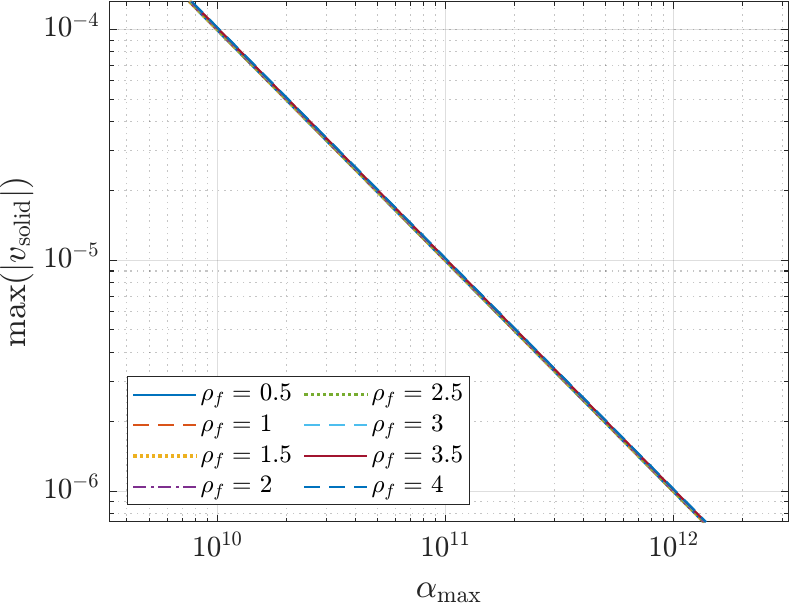}}
	\caption{Maximum velocity in the solid regions vs $\alpha_\text{max}$ for different $\rho_f$ values on a log-log scale.}
	\label{fig:kmax_vs_rho}
\end{figure*}

\subsection{Relation between $\alpha_\text{max}$ and $\rho_f$}
\label{ssec:alpha_rho_chrctr}

\begin{figure}
	\centering
	\includegraphics[width=0.5\textwidth]{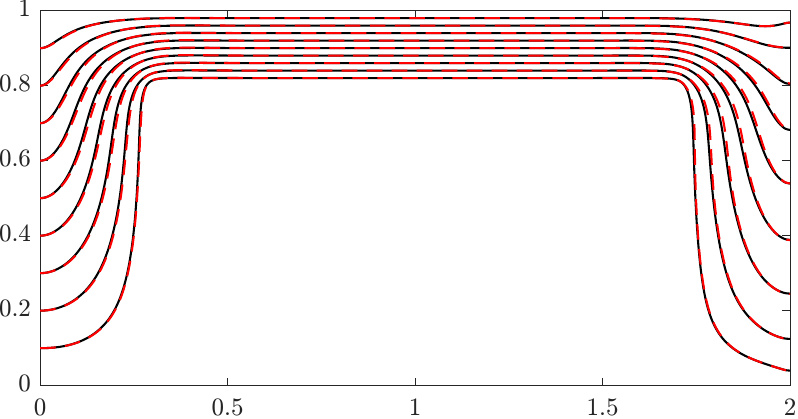}
	\caption{Comparison of velocity streamlines for $\rho_f = 0.5$ kg/m\textsuperscript{3} (solid black) vs $\rho_f = 4$ kg/m\textsuperscript{3} (dashed red).}
	\label{fig:v_strmlns_cmprsn}
\end{figure}

A fluid flow analysis is run for all combinations of the following values; $\rho_f = 0.5, 1, 1.5, ..., 4$ and $\alpha_\text{max} = 0, 1e0, 1e1,$ $..., 1e20$. The extracted results are presented in Fig. \ref{fig:kmax_vs_rho}, where it's clear that $\alpha_\text{max}$ is ``almost" independent of $\rho_f$. In fact, $\alpha_\text{max}$ is not entirely independent of $\rho_f$ due to the appearance of velocity components in the convective term on the left hand side of Eq. \ref{eq:nondmnsl_ns}. From a physics perspective, altering the value of $\rho_f$ affects the velocity field due to the changing ratio of inertia vs viscous forces. In Fig. \ref{fig:v_strmlns_cmprsn}, a comparison is presented between the velocity streamlines for the cases of $\rho_f = 0.5$ kg/m\textsuperscript{3} vs $\rho_f = 4$ kg/m\textsuperscript{3}, where the later case shows a slightly increased fluid inertia. Nonetheless, this effect is minimal on the velocity in the solid regions due to the fact that, in those regions, the viscous forces (i.e. contributions from Brinkman penalization and fluid viscosity) are much larger than the inertia forces. Hence, unless the change in $\rho_f$ exceeds an order of magnitude or Reynolds number is generally large, it's safe to ignore its effect on $\alpha_\text{max}$ from a practical perspective.

\subsection{Relation between $\alpha_\text{max}$ and $\mu$}
\label{ssec:kmax_mu_chrctr}

\begin{figure*}
	\centering
	\subfloat[]{\includegraphics[width=0.5\textwidth]{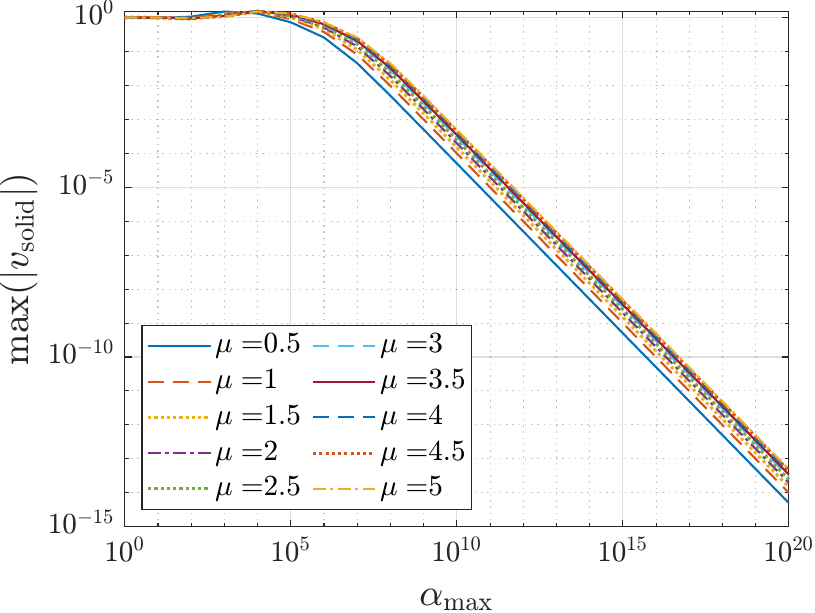}}
	\subfloat[]{\includegraphics[width=0.5\textwidth]{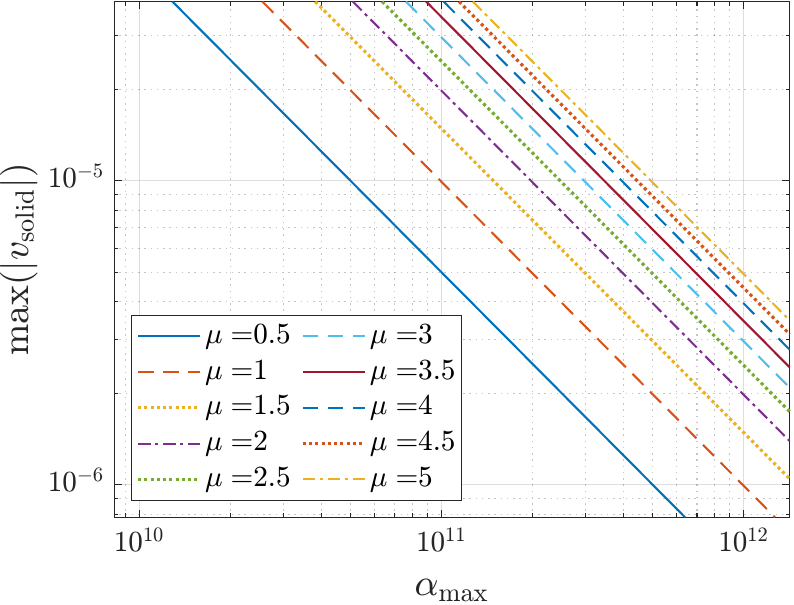}}
	\caption{Maximum velocity in the solid regions vs $\alpha_\text{max}$ for different $\mu$ values on a log-log scale.}
	\label{fig:kmax_vs_mu}
\end{figure*}

\begin{figure*}[t!]
	\captionsetup{margin={10pt}}
	\centering
	\subfloat[A log-log scale is used. From bottom up, max$(|v_\mathrm{solid}|) = 1e-2, 1e-3, ..., 1e-12$.]{\includegraphics[width=0.5\textwidth]{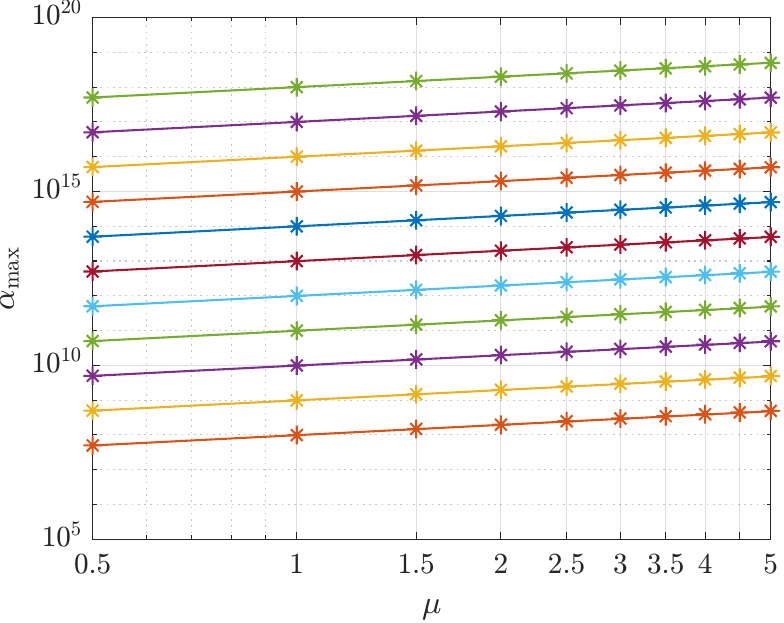}}
	\subfloat[A single case of max$(|v_\mathrm{solid}|) = 1e-12$ with error bars calculated as abs($\alpha_\text{max}|$\textsubscript{Data} - $\alpha_\text{max}|$\textsubscript{Eq.}).] {\includegraphics[width=0.5\textwidth]{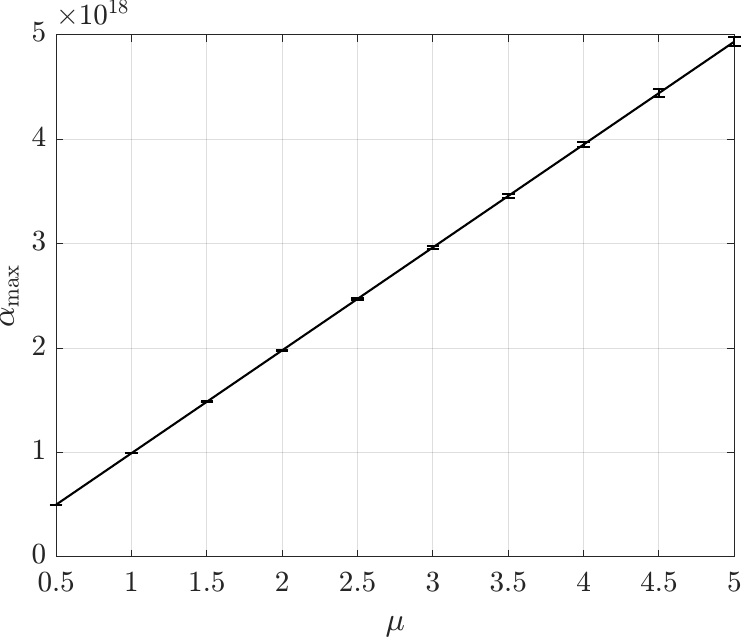}}
	\caption{Comparison of Eq. \ref{eq:kmax_vs_mu_crv_ftng} (asterisks) to data points from the numerical experiments (solid lines).}
	\label{fig:kmax_vs_mu_crv_ftng_chck}
\end{figure*}

A fluid flow analysis is run for all combinations of the following values; $\mu = 0.5, 1.0, 1.5, ..., 5.0$ and $\alpha_\text{max} = 0, 1e0, 1e1,$ $..., 1e20$. The extracted results are presented in Fig. \ref{fig:kmax_vs_mu}. Similarly to the approach followed in Section \ref{ssec:kmax_h_chrctr}, we limit the data used for curve fitting to only 4 points (2 points along $\mu$ and 2 points along max$(|v_\mathrm{solid}|)$); which are all combinations of $\mu = 0.5, 5.0$ and $\alpha_\text{max} = 1e8, 1e20$. The following relation is obtained:
\begin{equation}
	\alpha_\text{max} = 10^{-q}  \left( 9.857e5 \ \mu + 7331 \right)
	\label{eq:kmax_vs_mu_crv_ftng}
\end{equation}

\noindent where $q$ is defined similarly to Section \ref{ssec:kmax_h_chrctr}. To check the soundness of this relation, we compare Eq. \ref{eq:kmax_vs_mu_crv_ftng} to the original set of data points for $\mu = 0.5, 1.0, 1.5, ..., 5.0$ and $\alpha_\text{max} = 1e8, 1e9, 1e10, $ $..., 1e20$. The comparison is presented in Fig. \ref{fig:kmax_vs_mu_crv_ftng_chck} showing good agreement. Noticing that the error is consistently increasing with increasing $\mu$, we conjecture this error is due to the changing ratio in inertia vs viscous forces discussed in Section \ref{ssec:alpha_rho_chrctr}. Nonetheless, the maximum error at $\mu = 5$ Pa$\cdot$s is less than 1\% in the case of max$(|v_\mathrm{solid}|) = 10^{-12}$.

\subsection{Relation between $\alpha_\text{max}$ and $L_c$}
\label{ssec:kmax_lc_chrctr}

\begin{figure*}[t!]
	\centering
	\subfloat[]{\includegraphics[width=0.5\textwidth]{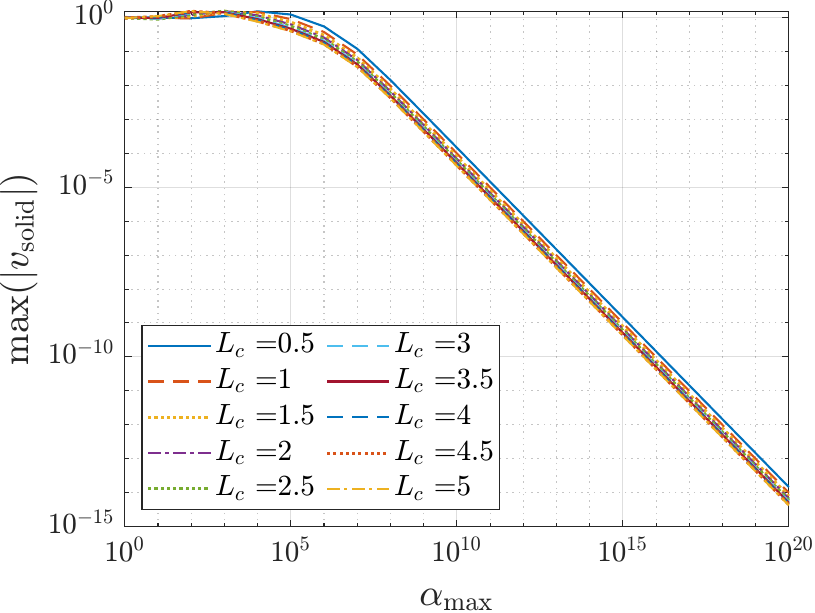}}
	\subfloat[]{\includegraphics[width=0.5\textwidth]{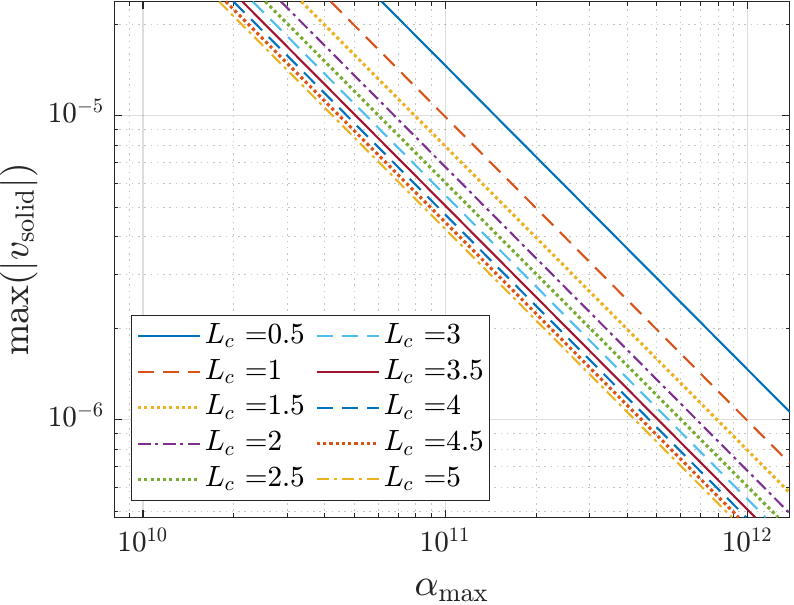}}
	\caption{Maximum velocity in the solid regions vs $\alpha_\text{max}$ for different $L_c$ values on a log-log scale.}
	\label{fig:kmax_vs_lc}
\end{figure*}

\begin{figure*}[t!]
	\captionsetup{margin={10pt}}
	\centering
	\subfloat[A log-log scale is used. From bottom up, max$(|v_\mathrm{solid}|) = 1e-2, 1e-3, ..., 1e-12$.]{\includegraphics[width=0.5\textwidth]{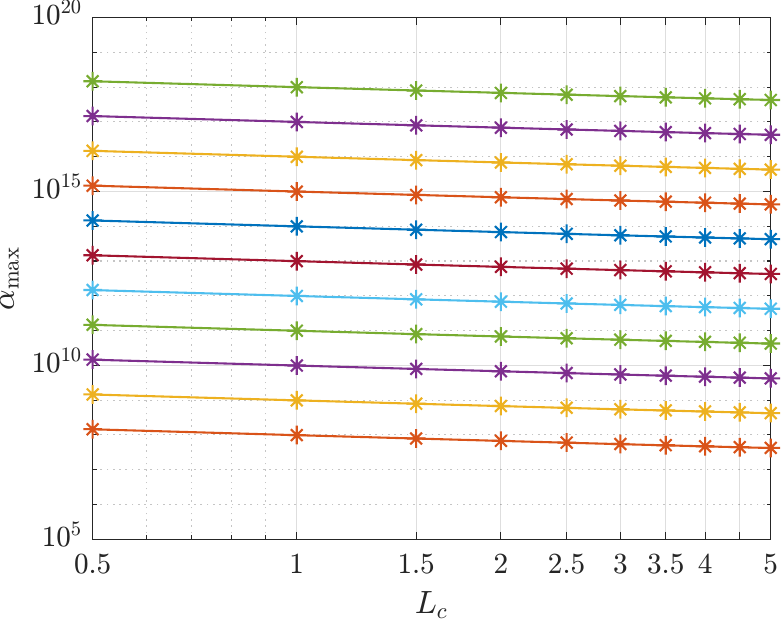}}
	\subfloat[A single case of max$(|v_\mathrm{solid}|) = 1e-12$ with error bars calculated as abs($\alpha_\text{max}|$\textsubscript{Data} - $\alpha_\text{max}|$\textsubscript{Eq.}).] {\includegraphics[width=0.5\textwidth]{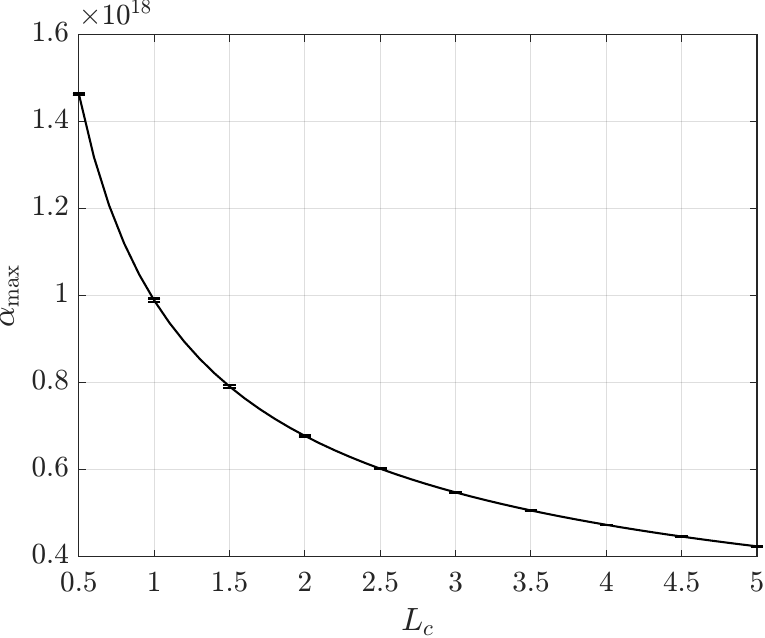}}
	\caption{Comparison of Eq. \ref{eq:kmax_vs_lc_crv_ftng} (asterisks) to data points from the numerical experiments (solid lines).}
	\label{fig:kmax_vs_lc_crv_ftng_chck}
\end{figure*}

A study is run for all combinations of the following values; $L_c = 0.5, 1.0, 1.5, ..., 5.0$ and $\alpha_\text{max} = 0, 1e0, 1e1,$ $..., 1e20$. The extracted results are presented in Fig. \ref{fig:kmax_vs_lc}. At first, we attempted to follow an approach similar to the one followed in Section \ref{ssec:kmax_h_chrctr} by limiting the data used for curve fitting to only 6 points (3 points along $L_c$ and 2 points along max$(|v_\mathrm{solid}|)$); which are all combinations of $L_c = 0.5, 3.0, 5.0$ and $\alpha_\text{max} = 1e8, 1e20$. However, the fitted equation in the form of $\alpha_\text{max} \propto 1/L_c^2$ showed considerable disagreement with the original data set extracted from numerical experiments. Secondly, we even attempted to use all data points in the curve fitting process, but still failed to get a satisfying agreement. This issue led us to believe that the use of $L_c$ in non-dimensionalizing $\alpha(\rho)$ in Eq. \ref{eq:nondmnsl_alpha} is incorrect. To gain some insight into the relation between $\alpha_\text{max}$ and $L_c$, we attempted to fit an equation in the form of $\alpha_\text{max} \propto a_1/{L_c}^{a_2}$ where $a_1$ and $a_2$ are constants. The following relation is obtained:
\begin{equation}
	\alpha_\text{max} = 10^{-q} \left( \frac{9.065e5}{{L_c}^{0.6073}} + 8.3e4 \right)
	\label{eq:kmax_vs_lc_crv_ftng}
\end{equation}

\noindent where $q$ is defined in Section \ref{ssec:kmax_h_chrctr}. Notice that in fitting this relation, we only used 6 points as discussed earlier. To check the soundness of this relation,  we compare Eq. \ref{eq:kmax_vs_lc_crv_ftng} to the original set of data points for $L_c = 0.5, 1.0, 1.5, ..., 5.0$ and $\alpha_\text{max} = 1e8, 1e9,$ $..., 1e20$. The comparison is presented in Fig. \ref{fig:kmax_vs_lc_crv_ftng_chck} showing surprisingly good agreement with a maximum error less than 0.5\% in the case of max$(|v_\mathrm{solid}|) = 10^{-12}$. We conjecture that $\alpha_\text{max}$ is in fact related to a length characteristic of the porous microstructure, and this length is related to $L_c$ on the macroscale.

\subsection{Relation between $\alpha_\text{max}$ and $v_c$}
\label{ssec:kmax_vc_chrctr}

\begin{figure*}
	\centering
	\subfloat[]{\includegraphics[width=0.5\textwidth]{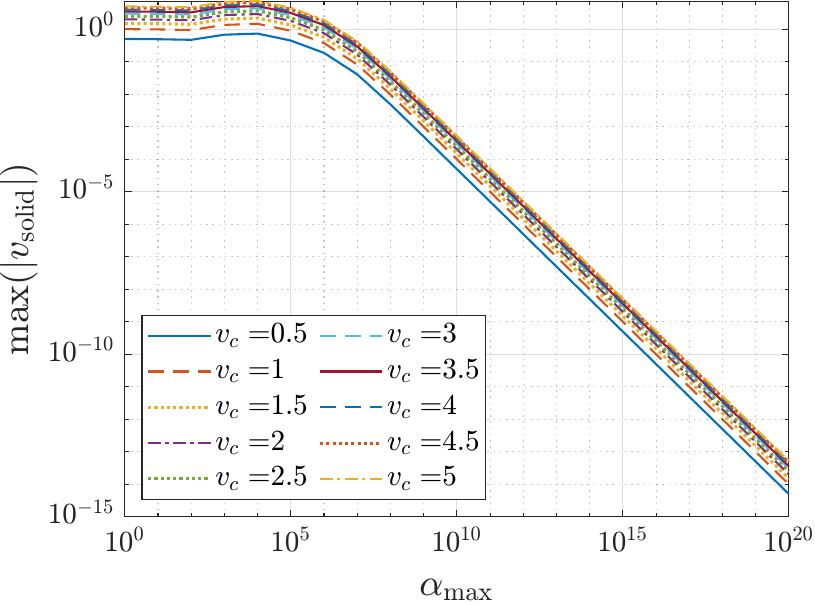}}
	\subfloat[]{\includegraphics[width=0.5\textwidth]{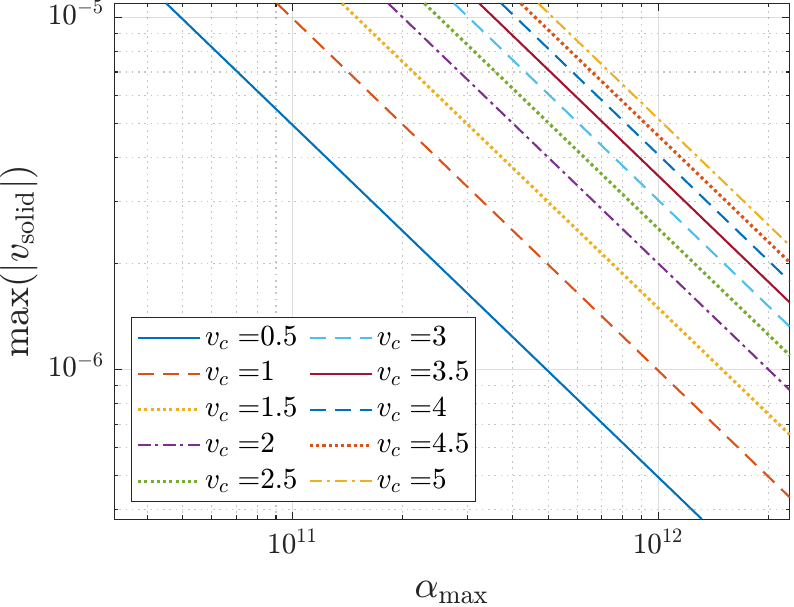}}
	\caption{Maximum velocity in the solid regions vs $\alpha_\text{max}$ for different $v_c$ values on a log-log scale.}
	\label{fig:kmax_vs_vc}
\end{figure*}

\begin{figure*}[t!]
	\captionsetup{margin={10pt}}
	\centering
	\subfloat[A log-log scale is used. From bottom up, max$(|v_\mathrm{solid}|) = 1e-2, 1e-3, ..., 1e-12$.]{\includegraphics[width=0.5\textwidth]{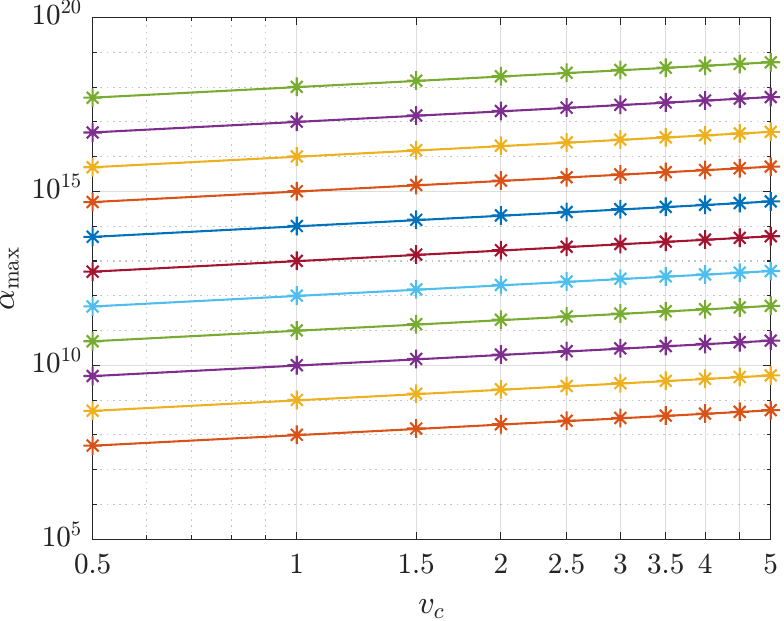}}
	\subfloat[A single case of max$(|v_\mathrm{solid}|) = 1e-12$ with error bars calculated as abs($\alpha_\text{max}|$\textsubscript{Data} - $\alpha_\text{max}|$\textsubscript{Eq.}).] {\includegraphics[width=0.5\textwidth]{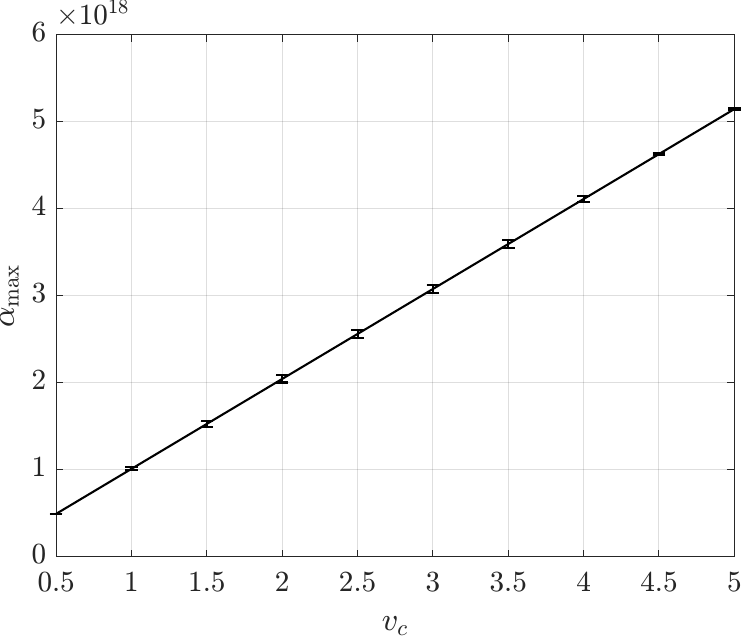}}
	\caption{Comparison of Eq. \ref{eq:kmax_vs_vc_crv_ftng} (asterisks) to data points from the numerical experiments (solid lines).}
	\label{fig:kmax_vs_vc_crv_ftng_chck}
\end{figure*}

A study is run for all combinations of the following values; $v_c = 0.5, 1.0, 1.5, ..., 5.0$ and $\alpha_\text{max} = 0, 1e0, 1e1,$ $..., 1e20$. The extracted results are presented in Fig. \ref{fig:kmax_vs_vc}. Similarly to the approach followed in Section \ref{ssec:kmax_h_chrctr}, we limit the data used for curve fitting to only 4 points (2 points along $v_c$ and 2 points along max$(|v_\mathrm{solid}|)$); which are all combinations of $v_c = 0.5, 5.0$ and $\alpha_\text{max} = 1e8, 1e20$. The following relation is obtained:
\begin{equation}
	\alpha_\text{max} = 10^{-q} \left(1.034e6 \ v_c - 2.253e4 \right)
	\label{eq:kmax_vs_vc_crv_ftng}
\end{equation}

\noindent where $q$ is defined similarly to Section \ref{ssec:kmax_h_chrctr}. To check the soundness of this relation, we compare Eq. \ref{eq:kmax_vs_vc_crv_ftng} to the original set of data points for $v_c = 0.5, 1.0, 1.5, ..., 5.0$ and $\alpha_\text{max} = 1e8, 1e9, 1e10, $ $..., 1e20$. The comparison is presented in Fig. \ref{fig:kmax_vs_vc_crv_ftng_chck} showing good agreement. Similarly to Section \ref{ssec:kmax_mu_chrctr}, the error appears to be decreasing with increasing $v_c$, we conjecture this error is due to the changing ratio in inertia vs viscous forces discussed in Section \ref{ssec:alpha_rho_chrctr}. Nonetheless, the maximum error is only 2.1\% in the case of max$(|v_\mathrm{solid}|) = 10^{-12}$.

\section{Conclusions}
\label{sec:conc}

In this work, we investigated the dependency of the inverse permeability maximum limit on the mesh size and flow conditions. The motivation behind this study is the need for mimicking the same behavior of the Brinkman-penalized Navier-Stokes equations for different mesh sizes and flow conditions, which is particularly useful when calibrating the various interpolation and projection parameters common in density-based topology optimization of fluid-dependent problems.

We first started by investigating the fluid flow governing equations in their strong as well as discretized finite element forms. We analytically derived proportionality relations between the maximum inverse permeability limit and the mesh size and flow condition parameters. We emphasize that these proportionality relations are not closed-form, instead they are generally true with a certain range of flow behavior. In general, these proportionality relations are independent of the design problem, though the proportionality coefficients are problem-dependent. 

For a specific design problem common in topology optimization of fluid-structure interactions, we proved these dependency relations numerically for the mesh size, dynamic viscosity, and characteristic velocity. For the characteristic length, a different relation was obtained from curve fitting which we believe is due to the dependency of the maximum inverse permeability limit on a microscale characteristic length that is somehow related to the macroscale one. In the case of the fluid density, it is deduced analytically and proven numerically that the maximum inverse permeability limit is independent of the fluid density when the change is within a reasonable range of Reynolds numbers. We also showed that only a handful of data points are needed to calculate proportionality coefficients for other problems, given that the analytical dependency relations are known a priori.

\backmatter

%\bmhead{Supplementary information}
%
%If your article has accompanying supplementary file/s please state so here. 
%
%Authors reporting data from electrophoretic gels and blots should supply the full unprocessed scans for key as part of their Supplementary information. This may be requested by the editorial team/s if it is missing.
%
%Please refer to Journal-level guidance for any specific requirements.
%
%\bmhead{Acknowledgements}
%
%Acknowledgements are not compulsory. Where included they should be brief. Grant or contribution numbers may be acknowledged.
%
%Please refer to Journal-level guidance for any specific requirements.

\section*{Declarations}

\noindent \textbf{Author Contributions} MA co-conceptualized the study, co-derived the analytical formulations, produced the results and graphs, wrote the draft, and co-wrote the conclusions. AC co-conceptualized the study, co-derived the analytical formulations, supervised, provided resources, managed and coordinated the research effort, reviewed and edited the draft, and co-wrote the conclusions.
\vspace{1em}

\noindent \textbf{Data Availability and Replication of Results} The authors included all data and described the formulations and algorithms needed to replicate the results. The authors will be happy to provide any additional data upon request.
\vspace{1em}

\noindent \textbf{Conflict of Interest} The authors declare they have no conflict of interest.
\vspace{1em}

\noindent \textbf{Funding} There is no funding source.
\vspace{1em}

\noindent \textbf{Ethical Approval} This article does not contain any studies with human participants or animals performed by any of the authors.

\bibliographystyle{apalike}
\bibliography{library}% common bib file
%% if required, the content of .bbl file can be included here once bbl is generated
%%\input sn-article.bbl

\end{document}